\documentclass[12pt]{article}%
\usepackage{amsfonts}
\usepackage{amssymb}
\usepackage{amsmath}%
\usepackage{graphicx}
\usepackage{hyperref}

\providecommand{\U}[1]{\protect\rule{.1in}{.1in}}

\newtheorem{theorem}{Theorem}

\newtheorem{corollary}[theorem]{Corollary}

\begin{document}

\title{Schur's exponent conjecture II}
\author{Michael Vaughan-Lee}
\date{December 2021}
\maketitle

\begin{abstract}
Primoz Moravec published a very important paper in 2007 where he proved that
if $G$ is a finite group of exponent $n$ then the exponent of the Schur
multiplier of $G$ can be bounded by a function $f(n)$ depending only on $n$.
Moravec does not give a value for $f(n)$, but actually his proof shows that we
can take $f(n)=ne$ where $e$ is the order of $b^{-n}a^{-n}(ab)^{n}$ in the Schur
multiplier of $R(2,n)$. (Here $R(2,n)$ is the largest finite two generator
group of exponent $n$, and we take $a,b$ to be the generators of $R(2,n)$.) It
is an easy hand calculation to show that $e=n$ for $n=2,3$, and it is a
straightforward computation with the $p$-quotient algorithm to show that $e=n$
for $n=4,5,7$. The groups $R(2,8)$ and $R(2,9)$ are way out of range of the
$p$-quotient algorithm, even with a modern supercomputer. But we are able to
show that $e\geq n$ for $n=8,9$. Moravec's proof also shows that if $G$ is a
finite group of exponent $n$ with nilpotency class $c$, then the exponent of
the Schur multiplier of $G$ is bounded by $ne$ where $e$ is the order of
$b^{-n}a^{-n}(ab)^{n}$ in the Schur multiplier of the class $c$ quotient
$R(2,n;c)$ of $R(2,n)$. If $q$ is a prime power we let $e_{q,c}$ be the order
of $b^{-q}a^{-q}(ab)^{q}$ in the Schur multiplier of $R(2,q;c)$. We are able
to show that $e_{p^{k},p^{2}-p-1}$ divides $p$ for all prime powers $p^{k}$.
If $k>2$ then $e_{2^{k},c}$ equals 2 for $c<4$, equals 4 for $4\leq c\leq11$,
and equals $8$ for $c=12$. If $k>1$ then $e_{3^{k},c}$ equals 1 for $c<3$,
equals 3 for $3\leq c<12$, and equals 9 for $c=12$.

We also investigate the order of $[b,a]$ in a Schur cover for $R(2,q;c)$.

\end{abstract}

\section{Introduction}

There is a long-standing conjecture attributed to I. Schur that the exponent
of the Schur multiplier, $M(G)$, of a finite group $G$ divides the exponent of
$G$. It is easy to show that this conjecture holds true for groups of exponent
2 and exponent 3, but a counterexample in exponent 4 was found by Bayes,
Kautsky and Wamsley \cite{bayes74} in 1974. The conjecture remained open for
odd exponent until 2020, when I found counterexamples of exponent 5 and
exponent 9 \cite{vleeschur}. It seems certain that there are counterexamples
to this conjecture for all prime powers greater than 3, but this leaves open
the question of what bounds on the exponent of $M(G)$ might hold true.

The usual definition of the Schur multiplier of a finite group $G$ is the
second homology group $H_{2}(G,\mathbb{Z})$. For computational purposes we use
the Hopf formulation of the Schur multiplier, which is as follows. We write
$G=F/R$, where $F$ is a free group, and then%
\[
M(G)=(R\cap F^{\prime})/[F,R].
\]
If we let $H=F/[F,R]$ then $H$ is an infinite group, and the quotient
$H/H^{\prime}$ is a free abelian group with rank equal to the rank of $F$ as a
free group. However the centre of $H$ contains $R/[F,R]$ and has finite index
in $H$. So the derived group $H^{\prime}=F^{\prime}/[F,R]$ is finite, and
$M(G)$ is finite. It is known that the exponent of $M(G)$ divides the order of
$G$. Furthermore if $G$ has exponent $n$ and if $h\in H^{\prime}$ then
$h^{n}\in M(G)$. So the quotient $H^{\prime}/M(G)$ has exponent dividing $n$.
In particular, if $G$ is a finite $p$-group then $H^{\prime}$ is a finite $p$-group.

Primoz Moravec published a very important paper \cite{moravec2007} in 2007 in
which he proved that if $G$ is a finite group of exponent $n$ then the
exponent of $M(G)$ is bounded by a function $f(n)$ depending only on $n$.
Moravec does not give an explicit formula for $f(n)$, but his proof of this
theorem actually shows that if $G$ is a finite group of exponent $n$ then the
exponent of $M(G)$ divides $ne$ where $e$ is the order of $b^{-n}%
a^{-n}(ab)^{n}$ in the Schur multiplier of $R(2,n)$. (We let $R(2,n)$ denote
the largest finite $2$ generator of exponent $n$, and we take the generators
of $R(2,n)$ to be $a,b$.)

\begin{theorem}
[Moravec, 2007]Let $b^{-n}a^{-n}(ab)^{n}\in M(R(2,n))$ have order $e$. If $G$
is any finite group of exponent $n$ then the exponent of $M(G)$ divides $ne$.
\end{theorem}

We give a short proof of this theorem in Section 2. The same proof gives the
following theorem.

\begin{theorem}
Let $R(2,n;c)$ be the nilpotent of class $c$ quotient of $R(2,n)$, and let
$e_{n,c}$ be the order of $b^{-n}a^{-n}(ab)^{n}\in M(R(2,n;c))$. If $G$ is any
finite group of exponent $n$ with class $c$, then the exponent of $M(G)$
divides $ne_{n,c}$.
\end{theorem}

Our proof of Moravec's theorem also gives the following corollary.

\begin{corollary}
Let $R(d,n)$ be the largest finite $d$ generator group of exponent $n$. Then
if $d\geq2$ the exponent of $M(R(d,n))$ is the order of $b^{-n}a^{-n}%
(ab)^{n}\in M(R(2,n))$.
\end{corollary}

Theorem 1 led me to investigate the order $e_{q}$ of $b^{-q}a^{-q}(ab)^{q}\in
M(R(2,q))$ for prime power exponents $q$. (We restrict ourselves to groups of
prime power exponent since if $G$ is any finite group and $p$ is any prime
then the Sylow $p$-subgroup of $M(G)$ is a subgroup of the Schur multiplier of
the Sylow $p$-subgroup of $G$.) It is an easy hand calculation to show that
$e_{q}=q$ for $q=2,3$. And it is a straightforward computation with the
$p$-quotient algorithm \cite{havasnew80} to show that $e_{q}=q$ for $q=4,5,7$.
Computing the groups $R(2,8)$ or $R(2,9)$ (or their Schur covers) is way out
of the range of the $p$-quotient algorithm, even with a modern supercomputer.
But I am able to show that $e_{q}\geq q$ for $q=8,9$.

Theorem 2 led me to investigate the order $e_{q,c}$ of $b^{-q}a^{-q}%
(ab)^{q}\in M(R(2,q;c))$ for prime power exponents $q$ and various $c$.

\begin{theorem}
If $q$ is a power of the prime $p$ then $e_{q,p^{2}-p-1}$ divides $p$.
\end{theorem}

\begin{theorem}
If $k>2$ then $e_{2^{k},c}$ equals $2$ for $c<4$, equals $4$ for $4\leq
c\leq11$, and equals $8$ for $c=12$.
\end{theorem}

\begin{theorem}
If $k>1$ then $e_{3^{k},c}$ equals $1$ for $c<3$, equals $3$ for $3\leq c<12$,
and equals $9$ for $c=12$.
\end{theorem}

The counterexamples to Schur's conjecture found in \cite{bayes74} and
\cite{vleeschur} are based on the following construction. Let $H(q,c)$ be the
largest four generator group of exponent $q$ and nilpotency class $c$ which is
generated by $a,b,c,d$ and subject to the relation%
\[
\lbrack b,a][d,c]=1.
\]
The Bayes, Kautsky and Wamsley example in \cite{bayes74} is $H(4,4)$ which has
Schur multiplier of exponent 8. The element $[b,a][d,c]$ has order 8 in the
Schur multiplier, and in fact $[b,a]$ has order $8$ in a Schur cover of $R(2,4;4)$.
Similarly my examples in \cite{vleeschur} are based on $H(5,9)$ and $H(9,9)$.
In the Schur multipliers of these two groups the elements $[b,a][d,c]$ have
order 25 and 27. The examples \textquotedblleft work\textquotedblright%
\ because $[b,a]$ has order 25 and 27 in Schur covers of $R(2,5;9)$
and $R(2,9;9)$. It seems plausible that more generally the exponent of
$M(H(q,c))$ is the order of $[b,a]$ in a Schur cover of $R(2,q;c)$, though I have no idea
how to prove this in general. For this reason I have looked at the order of
$[b,a]$ in Schur covers of $R(2,q;c)$ for various $q,c$. For example, $[b,a]$
has order 32 in a Schur cover of $R(2,8;12)$ so it seems to me to be 
extremely likely that $M(H(8,12))$ has exponent 32, though I cannot prove it
(yet!). This is an interesting example because all $p$-group counterexamples
$G$ to Schur's conjecture found so far have exp$\,M(G)=p\,$exp$\,G$.

I am able to prove that $M(H(q,4))$ has exponent $2q$ whenever $q\ge 4$
is a power of 2. For $q=9,27$ $M(H(q,9))$ has exponent $3q$, and it seems
very likely that $M(H(q,9))$ has exponent $3q$ for all $q>3$ which are powers
of 3. Similarly for $q=5,25$ $M(H(q,9))$ has exponent $5q$, and it seems
very likely that $M(H(q,9))$ has exponent $5q$ for all $q$ which are powers
of 5.

\begin{theorem}
If $q>4$ is a power of $2$, and if we let $f$ be the order of $[b,a]$ in a
Schur cover of $R(2,q;c)$ then $f=q$ for $c<4$, $f=2q$ for $4\leq c<12$, and
$f=4q$ for $c=12$. If $G$ is a finite $2$-group with nilpotency class less than
$4$ then exp$\,M(G)$ divides exp$\,G$.
\end{theorem}

\begin{theorem}
If $q>3$ is a power of $3$, and if we let $f$ be the order of $[b,a]$ in a
Schur cover of $R(2,q;c)$ then $f=q$ for $c<9$, and $f=3q$ for $c=9$.
If $G$ is any group of exponent $q$ with class less than $9$ then 
exp$\,M(G)$ divides $q$.
\end{theorem}

\begin{theorem}
Let $q$ be a power of $5$. If we let $f$ be the order of $[b,a]$ in a Schur
cover of $R(2,q;c)$ then $f=q$ for $c<9$, and $f=5q$ for $c=9$. If $G$ is any
group of exponent $q$ with class less than $9$ then exp$\,M(G)$ divides $q$.
\end{theorem}

\begin{theorem}
Let $q$ be a power of $7$. If we let $f$ be the order of $[b,a]$ in a Schur
cover of $R(2,q;c)$ then $f=q$ for $c<13$, and $f=7q$ for $c=13$. If $G$ is
any group of exponent $q$ with class less than $13$ then exp$\,M(G)$ divides
$q$.
\end{theorem}

Theorem 9 follows from the fact that in a Schur cover of $R(2,5)$ the element
$[b,a]^{5}$ is a product of commutators with entries $a$ or $b$, where at
least 5 of the entries are $a$'s and at least 5 of the entries are $b$'s.
Similarly Theorem 10 follows from the fact that in a Schur cover of $R(2,7)$
the element $[b,a]^{7}$ is a product of commutators with at least 7 entries
$a$ and at least 7 entries $b$. If the same pattern is repeated for higher
primes $p$, then I would expect the order of $[b,a]$ in $M(R(2,p^{k};c))$ to
be $p^{k}$ for $c<2p-1$ and to be $p^{k+1}$ when $c=2p-1$.

The theorems above omit the exponents 2,3,4. It is easy to see that the Schur
multiplier of a group of exponent 2 or exponent 3 has exponent 2 or 3
(respectively). Moravec \cite{moravec2007} proves that the Schur multiplier of
a group of exponent 4 has exponent dividing 8.

\section{Proof of Theorem 1}

We write $R(2,n)=F/R$ where $F$ is the free group generated by $a,b$, and let
$H=F/[F,R]$. (We are not assuming here that $n$ is a prime power.) Then
$b^{-n}a^{-n}(ab)^{n}\in M(R(2,n))$. Let $b^{-n}a^{-n}(ab)^{n}$ have order 
$e$. We show that $e$ is the exponent of $M(R(d,n))$ for all $d\geq2$, and
that if $G$ is any finite group of exponent $n$ then exp$(M(G))$ divides $ne$.

So let $G$ be \emph{any} finite group of exponent $n$, let $G=F/R$ where $F$
is a free group, and let $H=F/[F,R]$. (Apologies for using the same notation
for the covering group of $G$ as I used for the covering group of $R(2,n)$.)
Let $a,b$ be \emph{any} two elements in $H$. Then the subgroup of $H$
generated by $a$ and $b$ is a homomorphic image of the cover of $R(2,n)$, and
so $b^{-n}a^{-n}(ab)^{n}\in H^{\prime}$ lies in the centre of $H$ and has
order dividing $e$. Let $F$ be freely generated by the set $X$ and let $K $ be
the subgroup of $H$ generated by the elements $x^{n}[F,R]$ ($x\in X$). Then
$K$ is a free abelian group which intersects $H^{\prime}$ trivially. If $w$ is
an arbitrary element of $F$ we can write $w=x_{1}x_{2}\ldots x_{k}$ for some
$k$ and some $x_{1},x_{2},\ldots,x_{k}\in X\cup X^{-1}$. Letting $a=x_{1}$ and
$b=x_{2}x_{3}\ldots x_{k}$ we see that%
\[
w^{n}=(x_{1}x_{2}\ldots x_{k})^{n}=x_{1}^{n}(x_{2}\ldots x_{k})^{n}%
b^{-n}a^{-n}(ab)^{n}.
\]
Repeating this argument we see that%
\[
w^{n}=x_{1}^{n}x_{2}^{n}\ldots x_{k}^{n}c,
\]
where $c$ is a product of terms of the form $b^{-n}a^{-n}(ab)^{n}$ with
$a,b\in F$. So we see that if $h\in H$ then $h^{n}$ is a product of an element
in $K$ and an element in $H^{\prime}$ which lies in the centre of $H $ and has
order dividing $e$. It follows that any product of $n^{th}$ powers in $H$ can
be expressed in the same form. Since $K\cap H^{\prime}=\{1\}$ we see that this
implies that any product of $n^{th}$ powers in $H$ which lies in $H^{\prime}$
has order dividing $e$. So exp$(M(R(d,n)))=e\,$\ for all $d\geq2$. If $h\in
H^{\prime}$, then $h^{n}$ is an $n^{th}$ power which lies in
$H^{\prime}$, and so $h^{ne}=1$. So $H^{\prime}$ has exponent dividing $ne$,
and this implies that exp$(M(G))$ divides $ne$.

\section{Some commutator calculus}

Let $F$ be the free group of rank 2 generated by $a$ and $b$. If we are
working in the nilpotent quotient $F/\gamma_{k+1}(F)$ for some $k$ then we pick
a fixed ordered set of basic commutators of weight at most $k$. See [3,
Chapter 11]. The first few basic commutators in our sequence are%
\[
a,\,b,\,[b,a],\,[b,a,a],\,[b,a,b],\,[b,a,a,a],\,[b,a,a,b],\,[b,a,b,b],\,[b,a,a,[b,a]].
\]
If $c_{1},c_{2},\ldots,c_{m}$ is our list of basic commutators of weight at
most $k$ then every element of $F/\gamma_{k+1}(F)$ can be written uniquely in
the form%
\[
c_{1}^{n_{1}}c_{2}^{n_{2}}\ldots c_{m}^{n_{m}}\gamma_{k+1}(F)
\]
for some integers $n_{1},n_{2},\ldots,n_{m}$. From the theory of Hall
collection (see [3, Theorem 12.3.1]), if $n$ is any positive
integer then in $F/\gamma_{k+1}(F)$
\begin{equation}
(ab)^{n}=a^{n}b^{n}[b,a]^{\binom{n}{2}}[b,a,a]^{\binom{n}{3}}c_{5}^{n_{5}%
}\ldots c_{m}^{n_{m}}%
\end{equation}
where the exponents $n,\binom{n}{2},\binom{n}{3},n_{5},\ldots,n_{m}$ take a very special form. If
$c_{r}$ has weight $w$ then $n_{r}$ is an integral linear combination of the
binomial coefficients $n,\binom{n}{2},\binom{n}{3},\ldots,\binom{n}{w}$.
Furthermore the integer coefficients which arise in these linear combinations
are positive, and are independent of $n$.

The exponents $n,\binom{n}{2},\binom{n}{3},n_{5},\ldots,n_{m}$ which arise in
equation (1) are all polynomials in $n$ which take integer values when $n$ is
an integer. The formula%
\[
\binom{n}{r}=\frac{n(n-1)\ldots(n-r+1)}{r!}
\]
also makes sense when $n$ is negative. We let $P=\mathbb{Q}[t]$ be the ring of
polynomials in an indeterminate $t$ over the rationals $\mathbb{Q}$. An
\emph{integer-valued polynomial} is a polynomial $f(t)\in P$ which takes
integer values whenever $f(t)$ is evaluated at an integer $n$. The set of
integer-valued polynomials is a subring of $P$, and is a free abelian group
with basis%
\[
1,\,t,\,\frac{t(t-1)}{2!},\,\ldots,\frac{t(t-1)\ldots(t-d+1)}{d!},\,\ldots.
\]
The exponents $n_{i}$ which arise in equation (1) all take the form
$n_{i}=f_{i}(n)$ where $f_{i}(t)$ is an integer-valued polynomial of degree at
most wt$\,c_{i}$ which does not depend on $n$. The polynomials $f(t)\in P$
that arise in this way in equation (1) also satisfy $f(0)=0$.

We rewrite equation (1) in the following form%
\begin{equation}
(ab)^{n}=a^{n}b^{n}[b,a]^{\binom{n}{2}}[b,a,a]^{\binom{n}{3}}c_{5}^{f_{5}%
(n)}\ldots c_{m}^{f_{m}(n)}%
\end{equation}
where the integer-valued polynomials $f_{i}(t)$ are independent of $n$.
The key properties of these polynomials to keep in mind are that
$f_i(0)=0$ and $\deg f_i(t)\leq \text{wt}\,c_i$.

We use equation (2) to get an expansion of $[y^{n},x]$ for $x,y\in F$.
Equation (2) gives%
\[
y^{n}x=x(y[y,x])^{n}=xy^{n}[y,x]^{n}[y,x,y]^{\binom{n}{2}}(c_{4}\alpha
)^{f_{4}(n)}\ldots(c_{m}\alpha)^{f_{m}(n)}\text{ modulo }\gamma_{k+1}(F)
\]
where $\alpha$ is the endomorphism of $F$ mapping $a,b$ to $y,[y,x]$. This
equation gives%
\begin{equation}
\lbrack y^{n},x]=[y,x]^{n}[y,x,y]^{\binom{n}{2}}(c_{4}\alpha)^{f_{4}(n)}%
\ldots(c_{m}\alpha)^{f_{m}(n)}\text{ modulo }\gamma_{k+1}(F).
\end{equation}

\section{Proof of Theorem 4}

We want to prove that if $q$ a power of the prime $p$ then the order of
$b^{-q}a^{-q}(ab)^{q}$ in $M(R(2,q;p^{2}-p-1)$ divides $p$. The case $p=2$ is
covered by Theorem 5, and so we assume that $p>2$. We write $R(2,q;p^{2}%
-p-1)=F/R$, where $F$ is the free group with free generators $a,b$. Let
$H=F/[F,R]$. So $H$ is nilpotent of class at most $p^{2}-p$, and (as we noted
in the introduction) $H^{\prime}$ is a finite $p$-group. We let
$c_{1},c_{2},\ldots,c_{m}$ be our list of basic commutators of weight at most
$p^{2}-p$ as described in Section 3. Let $x,y\in H$ and set $n=q$ in equation
(3) from Section 3. Using the fact that $H$ is nilpotent of class $p^{2}-p$ we
obtain a relation%
\begin{equation}
\lbrack y,x]^{q}[y,x,y]^{\binom{q}{2}}[y,x,y,y]^{\binom{q}{3}}(c_{5}%
\alpha)^{f_{5}(q)}\ldots(c_{m}\alpha)^{f_{m}(q)}=1,
\end{equation}
where $\alpha$ is the homomorphism from $F$ to $H$ mapping $a,b$ to $y,[y,x]
$. Recall that if wt$\,c_{i}=w$ then $\deg f_{i}(t)\leq w$. Note that if
wt$\,c_{i}=w$ then $c_{i}\alpha$ is a commutator in $x,y$ with $w$ entries
$y$. Also note that wt$\,c_{i}<p^{2}$ so that $f_{i}(q)$ is divisible by
$\frac{q}{p}$ for all $i$. And if wt$\,c_{i}<p$ then $f_{i}(q)$ is divisible
by $q$.

Since $H$ is nilpotent of class at most $p^{2}-p$, if $y\in\gamma_{p-1}(H)$
then $c_{i}\alpha=1$ whenever wt$\,c_{i}\geq p$. So if $y\in\gamma_{p-1}(H)$
then relation (4) shows that $[y,x]^{q}$ is a product of $q^{th}$ powers of
commutators in $x,y$ with higher weight. First let $y\in\gamma_{p^{2}-p-1}%
(H)$. Then equation (4) gives $[y,x]^{q}=1$. Since elements $[y,x]$ of this
form generate $\gamma_{p^{2}-p}(H)$ this implies that $\gamma_{p^{2}-p}(H)$
has exponent $q$. Next let $y\in\gamma_{p^{2}-p-2}(H) $. Then equation (4)
gives $[y,x]^{q}\in\gamma_{p^{2}-p}(H)^{q}=\{1\}$. So we see that
$\gamma_{p^{2}-p-1}(H)$ has exponent $q$. We continue in this way,
successively proving that $\gamma_{p^{2}-p-2}(H)$, $\gamma_{p^{2}-p-3}(H)$,
\ldots\ have exponent $q$. Finally we let $y\in\gamma_{p-1}(H)$ and prove that
$\gamma_{p}(H)$ has exponent $q$.

Now set $n=pq$ in equation (3), and we obtain the relation$.$%
\[
\lbrack y,x]^{pq}[y,x,y]^{\binom{pq}{2}}[y,x,y,y]^{\binom{pq}{3}}(c_{5}%
\alpha)^{f_{5}(pq)}\ldots(c_{m}\alpha)^{f_{m}(pq)}=1
\]
where all the exponents $f_{i}(pq)$ are divisible by $q$, and where
commutators with less than $p$ entries $y$ have exponents divisible by $pq$.
Since all commutators in $H$ with weight at least $p$ have order dividing $q$,
this implies that $[x,y]^{pq}$ is a product of $(pq)^{th}$ powers of commutators
of higher weight in $x,y$. So $[x,y]^{pq}=1$, and hence all commutators in $H$
have order dividing $pq$. Since $H$ has class less than $p^{2}$ and $\gamma_{p}(H)$ 
has exponent $q$ this implies that $H^{\prime}$ has exponent dividing $pq$.

Finally consider $(b^{-q}a^{-q}(ab)^{q})^{p}=b^{-pq}a^{-pq}(ab)^{pq}$.
Equation (2) from Section 3 gives%
\[
(b^{-q}a^{-q}(ab)^{q})^{p}=[b,a]^{\binom{pq}{2}}[b,a,a]^{\binom{pq}{3}}%
c_{5}^{f_{5}(pq)}\ldots c_{m}^{f_{m}(pq)}.
\]
All the exponents in this product are divisible by $q$, and the exponents of
commutators of weight less than $p$ in the product are divisible by $pq$. So
$(b^{-q}a^{-q}(ab)^{q})^{p}=1$.

The same argument (though easier) shows that in $M(R(2,q;p-2)$ the element
$b^{-q}a^{-q}(ab)^{q}=1$.

\section{Proof of Theorem 5 and Theorem 7}

Let $q=2^{k}$ $(k>2)$, and let $e_{q,c}$ be the order of $b^{-q}%
a^{-q}(ab)^{q}$ in $M(R(2,q;c))$. We want to prove that $e_{q,c}$ equals $2$
for $c<4$, equals $4$ for $4\leq c\leq11$, and equals $8$ for $c=12$. We also
want to prove that if $f$ is the order of $[b,a]$ in a Schur cover of
$R(2,q;c)$ then $f=q$ for $c<4$, $f=2q$ for $4\leq c<12$, and $f=4q$ when
$c=12$.

Let $R(2,q;12)=F/R$ where $F$ is the free group of rank two generated by $a,b
$, and let $H=F/[F,R]$. Then $H$ is an infinite group, but the subgroup
$\langle a^{q},b^{q}\rangle\leq H$ is a central subgroup with trivial
intersection with $H^{\prime}$. We can factor this subgroup out, without
impacting $M(R(2,q;12))$, and we now have a finite $2$-group. We used the
$p$-quotient algorithm to compute this quotient for $q=8,16,32$. (These were
quite easy computations.) The computations showed that $e_{q,c}$ takes the
values given in Theorem 5 for $q=8,16,32$, and that $f$ takes the values given
in Theorem 7 for $q=8,16,32$. We show that the fact that Theorem 5 and Theorem
7 hold true for $q=2^{5}$ implies that they hold true for all exponents
$q=2^{k}$ with $k\geq5$.

We let $q=2^{k}$ where $k\geq5$, and let $c_{1},c_{2},\ldots,c_{m}$ be our
list of basic commutators of weight at most $13$ as described in Section 3. As
in the proof of Theorem 4 we let $x,y\in H$ and obtain a relation%
\begin{equation}
\lbrack y^q,x]=[y,x]^{q}[y,x,y]^{\binom{q}{2}}[y,x,y,y]^{\binom{q}{3}}(c_{5}%
\alpha)^{f_{5}(q)}\ldots(c_{m}\alpha)^{f_{m}(q)}
\end{equation}
where $\alpha$ is the homomorphism from $F$ to $H$ mapping $a,b$ to $y,[y,x]
$. If wt$\,c_{i}=w$ then $c_{i}\alpha$ is a commutator in $x$ and $y$, with
$w$ entries $y$, and $\deg f_{i}\leq w$. The binomial coefficients
$\binom{q}{2}$ and $\binom{q}{3}$ are both divisible by $\frac{q}{2}$, the
binomial coefficients $\binom{q}{d}$ for $d<8$ are divisible by
$\frac{q}{4}$, and the binomial coefficients $\binom{q}{d}$ for $d<16$
are divisible by $\frac{q}{8}$.

If we let $y\in\gamma_{7}(H)$ then $[y,x,y]\in\gamma_{15}(H)=\{1\}$, so we see
that $[y,x]^{q}=1$. Now $\gamma_{8}(H)$ is generated by elements $[y,x]$ with
$y\in\gamma_{7}(H)$, and $\gamma_{8}(H)$ is abelian. So $\gamma_{8}(H)$ has
exponent $q$.

Now let $y\in\gamma_{4}(H)$, and replace $q$ by $2q$ in equation (5). Using
the fact that $\gamma_{8}(H)$ has exponent $q$ we see that $[y,x]^{2q}=1$. So
$\gamma_{5}(H)$ is generated by elements of order $2q$. Furthermore
$\gamma_{5}(H)$ is nilpotent of class 2, and $\gamma_{5}(H)^{\prime}\leq
\gamma_{10}(H)$ has exponent $q$. So $\gamma_{5}(H)$ has exponent $2q$.

Next let $y\in H^{\prime}$, and replace $q$ by $4q$ in equation (5). We obtain
$[y,x]^{4q}=1$. So $\gamma_{3}(H)$ is generated by elements of order $4q$, and
using the fact that $\gamma_{5}(H)$ has exponent $2q$ and $\gamma_{8}(H)$ has
exponent $q$, we see that $\gamma_{3}(H)$ has exponent $4q$.

Finally replace $q$ by $8q$ in equation (5) and we obtain $[y,x]^{8q}=1$ for
all $x,y$. Using facts that $\gamma_{3}(H)$ has exponent $4q$,
$\gamma_{5}(H)$ has exponent $2q$, and $\gamma_{8}(H)$ has
exponent $q$, we see that $H^{\prime}$ has exponent $8q$.

Let $N$ be the normal subgroup
\[
\gamma_{2}(F)^{8q}\gamma_{3}%
(F)^{4q}\gamma_{5}(F)^{2q}\gamma_{8}(F)^{q}\gamma_{14}(F)<F.%
\]
Then $H=F/M$ where $M=\langle\lbrack y^{q},x]\,:\,x,y\in F\rangle N$.

Next we let $K$ be the normal subgroup 
\[
\gamma_{2}(F)^{q}\gamma_{3}%
(F)^{\frac{q}{2}}\gamma_{5}(F)^{\frac{q}{4}}\gamma_{8}(F)^{\frac{q}{8}}%
\gamma_{14}(F)<F.
\]
(The relevance of $K$ is that if $q\geq8$ is a power of 2
then $[y^{q},x]\in K$ for all $x,y\in F$.) We show that every element $k\in K$
can be written uniquely modulo $\gamma_{14}(F)$ in the form%
\begin{equation}
k=[b,a]^{qm_{3}}c_{4}^{\frac{q}{2}m_{4}}\ldots c_{8}^{\frac{q}{2}m_{8}}%
c_{9}^{\frac{q}{4}m_{9}}\ldots c_{41}^{\frac{q}{4}m_{41}}c_{42}^{\frac{q}%
{8}m_{42}}\ldots c_{1377}^{\frac{q}{8}m_{1377}},%
\end{equation}
where $m_{3},m_{4},\ldots,m_{1377}$ are integers. (The number of basic
commutators of weight at most 4 is 8, the number of weight at most 7 is 41,
and the number of weight at most 13 is 1377.) First let%
\[
K_{2}=\gamma_{3}(F)^{\frac{q}{2}}\gamma_{5}(F)^{\frac{q}{4}}\gamma
_{8}(F)^{\frac{q}{8}}\gamma_{14}(F).
\]
We show that every element in $\gamma_{2}(F)^{q}$ can be written as
$[b,a]^{qm_{3}}$ modulo $K_{2}$. The elements of $\gamma_{2}(F)^{q}$ are
products of $q^{th}$ powers of elements in $\gamma_{2}(F)$. Let $x,y\in
\gamma_{2}(F)$. Then from equation (2) in Section 3 we see that%
\[
(xy)^{q}=x^{q}y^{q}[y,x]^{\binom{q}{2}}[y,x,x]^{\binom{q}{3}}(c_{5}%
\alpha)^{f_{5}(q)}\ldots(c_{1377}\alpha)^{f_{1377}(q)}
\]
where $\alpha:F\rightarrow F$ is the endomorphism mapping $a,b$ to $x,y$. Now
$[y,x]$, and $[y,x,x]\in\gamma_{4}(F)$ and $\binom{q}{2}$ and $\binom{q}{3}$
are divisible by $\frac{q}{2}$. So $[y,x]^{\binom{q}{2}}$ and $[y,x,x]^{\binom
{q}{3}}\in K_{2}$. Similarly all the terms $(c_{5}\alpha)^{f_{5}(q)}$,
$\ldots$, $(c_{1377}\alpha)^{f_{1377}(q)}$ lie in $K_{2}$. So $(xy)^{q}=x^{q}y^{q}$
modulo $K_{2}$. This means that every product of $q^{th}$ powers of elements
in $\gamma_{2}(F)$ can be written as the $q^{th}$ power of a single element in $\gamma_{2}(F)$
modulo $K_{2}$. So consider $x^{q}$ when $x\in\gamma_{2}(F)$. We can write
$x=[b,a]^{m_{3}}g$ for some $g\in\gamma_{3}(F)$. Then, as we have just seen,
$x^{q}=[b,a]^{qm_{3}}g^{q}$ modulo $K_{2}$, and since $g^{q}\in K_{2}$ this
means that $x^{q}=[b,a]^{qm_{3}}$ modulo $K_{2}$.

Now let $K_{3}=\gamma_{5}(F)^{\frac{q}{4}}\gamma_{8}(F)^{\frac{q}{8}}%
\gamma_{14}(F)$. Then $K_{2}$ is generated modulo $K_{3}$ by $(\frac{q}%
{2})^{th}$ powers of elements in $\gamma_{3}(F)$. Using the same argument as
above we see that if $x,y\in\gamma_{3}(F)$ then $(xy)^{\frac{q}{2}}%
=x^{\frac{q}{2}}y^{\frac{q}{2}}$ modulo $K_{3}$. So every product of $(\frac
{q}{2})^{th}$ powers in $\gamma_{3}(F)$ can be expressed modulo $K_{3}$ as
$x^{\frac{q}{2}}$ with $x\in\gamma_{3}(F)$. Let $x=c_{4}^{m_{4}}c_{5}^{m_{5}%
}\ldots c_{8}^{m_{8}}g$, with $g\in\gamma_{5}(F)$. Then%
\begin{align*}
x^{\frac{q}{2}}  & =c_{4}^{\frac{q}{2}m_{4}}(c_{5}^{m_{5}}\ldots c_{8}^{m_{8}%
}g)^{\frac{q}{2}}\text{ modulo }K_{3}\\
& =c_{4}^{\frac{q}{2}m_{4}}c_{5}^{\frac{q}{2}m_{5}}(c_{6}^{m_{6}}\ldots
c_{8}^{m_{8}}g)^{\frac{q}{2}}\text{ modulo }K_{3}\\
& =\ldots\\
& =c_{4}^{\frac{q}{2}m_{4}}\ldots c_{8}^{\frac{q}{2}m_{8}}g^{\frac{q}{2}%
}\text{ modulo }K_{3}\\
& =c_{4}^{\frac{q}{2}m_{4}}\ldots c_{8}^{\frac{q}{2}m_{8}}\text{ modulo }%
K_{3}.
\end{align*}
So every element of $\gamma_{2}(F)^{q}\gamma_{3}(F)^{\frac{q}{2}}$ can be
expressed modulo $K_{3}$ in the form%
\[
\lbrack b,a]^{qm_{3}}c_{4}^{\frac{q}{2}m_{4}}\ldots c_{8}^{\frac{q}{2}m_{8}}.
\]

Continuing in this way we see that every element $k\in K$ can be written in
the form (6) modulo $\gamma_{14}(F)$. Since every element of $F/\gamma
_{14}(F)$ can be uniquely expressed in the form%
\[
c_{1}^{n_{1}}c_{2}^{n_{2}}\ldots c_{1377}^{n_{1377}}\gamma_{14}(F)
\]
for some integers $n_{1},n_{2},\ldots,n_{1377}$, the expression (6) for $k$ modulo
$\gamma_{14}(F)$ is unique.

Similarly every element of $N$ can be expressed uniquely modulo $\gamma
_{14}(F)$ in the form%
\[
\lbrack b,a]^{8qn_{3}}c_{4}^{4qn_{4}}\ldots c_{8}^{4qn_{8}}c_{9}^{2qn_{9}%
}\ldots c_{41}^{2qn_{41}}c_{42}^{qn_{42}}\ldots c_{1377}^{qn_{1377}},
\]
and if $k\in K$ is given by (6) then $k\in N$ if and only if $8|m_{i}$ for
$3\leq i\leq1377$. So $K$ is generated by 
\[
C=\{[b,a]^{q},c_{4}^{\frac{q}{2}},\ldots,c_{8}^{\frac{q}{2}},c_{9}^{\frac{q}{4}},%
\ldots,c_{41}^{\frac{q}{4}},c_{42}^{\frac{q}{8}}, \ldots, %
c_{1377}^{\frac{q}{8}}\}
\]
modulo $\gamma_{14}(F)$, and all these generators have order 8 modulo $N$. We show
that provided $q\geq 32$ these generators commute with each other modulo $N$, so that $K/N$ is
abelian. It follows that $K/N$ is a direct sum of 1375 copies of the cyclic
group of order 8. If $k\in K$ has the form (6) then we let $[m_{3}%
,m_{4},\ldots,m_{1377}]$ be the \emph{representative vector} of $kN$, and we
think of this vector as an element in $C_{8}^{1375}$. Multiplying two elements
of $K/N$ corresponds to adding their representative vectors, and $k\in N$ if
and only if the representative vector of $kN$ is zero.

To show that the elements in $C$ commute with each other we first let $c,d\in\gamma_{3}(F)$ 
and consider the commutator $[c^r,d^s]$ for general $r,s>0$. The subgroup $\langle c,d\rangle$ 
has class at most 4 modulo $N$, and we expand $[c^{r},d^{s}]$ modulo $\gamma_{5}(\langle
c,d\rangle)$. Taking $x=d^{s}$ we see that
\[
\lbrack c^{r},x]=[c,x]^{r}[c,x,c]^{\binom{r}{2}}[c,x,c,c]^{\binom{r}{3}}.
\]
Also%
\[
\lbrack c,x]=[c,d^{s}]=[c,d]^{s}[c,d,d]^{\binom{s}{2}}[c,d,d,d]^{\binom{s}{3}%
}.
\]
So, modulo $\gamma_{5}(\langle c,d\rangle)$,
\begin{align*}
\lbrack c^{r},d^{s}]  & =([c,d]^{s}[c,d,d]^{\binom{s}{2}}[c,d,d,d]^{\binom
{s}{3}})^{r}([c,d]^{s}[c,d,d]^{\binom{s}{2}},c])^{\binom{r}{2}}[[c,d]^{s}%
,c,c]^{\binom{r}{3}}\\
& =[c,d]^{rs}[c,d,d]^{r\binom{s}{2}}[c,d,d,d]^{r\binom{s}{3}}[c,d,c]^{\binom
{r}{2}s}[c,d,d,c]^{\binom{r}{2}\binom{s}{2}}[c,d,c,c]^{\binom{r}{3}s}.
\end{align*}
Now assume that $q\ge 32$ and let $r=s=\frac{q}{2}$ then $2q$ divides $rs$ and all the other exponents
in the product above are divisible by $q$ so that $[c^{\frac{q}{2}}%
,d^{\frac{q}{2}}]\in N$. So the elements in $\{c_{i}^{\frac{q}{2}}%
:\,$wt$\,c_{i}=3,4\}$ commute with each other modulo $N$. Similarly if
$c\in\gamma_{3}(F)$ and $d\in\gamma_{5}(F)$ then $\langle c,d\rangle$ is
nilpotent of class at most 3 modulo $N$, and in the expansion of
\thinspace$\lbrack c^{\frac{q}{2}},d^{\frac{q}{4}}]$ the exponent of $[c,d]$
is divisible by $2q$, and all the exponents of terms of weight 3 are divisible
by $q$, so that $[c^{\frac{q}{2}},d^{\frac{q}{4}}]\in N$. So elements in
$\{c_{i}^{\frac{q}{2}}:\,$wt$\,c_{i}=3,4\}$ commute with elements in
$\{c_{i}^{\frac{q}{4}}\,:\,5\leq\,$wt$\,c_{i}\leq7\}$ modulo $N$. If
$c\in\gamma_{3}(F)$ and $d\in\gamma_{8}(F)$ then $\langle c,d\rangle$ is
nilpotent of class at most 2 modulo $N$, so that%
\[
\lbrack c^{\frac{q}{2}},d^{\frac{q}{8}}]=[c,d]^{\frac{q^{2}}{16}}\in N.
\]
In the same way, if $c,d\in\gamma_{5}(F)$ then $\langle c,d\rangle$ is
nilpotent of class at most 2 modulo $N$, and%
\begin{align*}
\lbrack c^{\frac{q}{4}},d^{\frac{q}{4}}]  & =[c,d]^{\frac{q^{2}}{16}}\in N,\\
\lbrack c^{\frac{q}{4}},d^{\frac{q}{8}}]  & =[c,d]^{\frac{q^{2}}{32}}\in N.
\end{align*}
Finally, if $c,d\in\gamma_{7}(F)$ then $[c,d]\in N$, so $[c^{r},c^{s}]\in N$
for all $r,s$. So all the elements in $C\backslash\{[b,a]^{q}\}$ commute with
each other modulo $N$ (provided $q\ge 32$).

It remains to show that $[b,a]^{q}$ commutes with elements in $C$. Let
$d\in\gamma_{3}(F)$ and let $c=[b,a]$. We obtain an expression for
$[d^{r},c^{s}]$ modulo $N$ similar to the expression we obtained above for
$[c^{r},d^{s}]$ modulo $N$ when $c,d\in\gamma_{3}(F)$. In this case we have
$\gamma_{7}(\langle c,d\rangle)\leq N$, so we obtain an expression involving
basic commutators in $c,d$ of weight at most 6. A complete list of basic
commutators up to weight 5 is%
\begin{align*}
& c,d,[d,c],[d,c,c],[d,c,d],[d,c,c,c],[d,c,c,d],[d,c,d,d],[d,c,c,c,c],\\
& [d,c,c,c,d],[d,c,c,d,d],[d,c,d,d,d],[d,c,c,[d,c]],[d,c,d,[d,c]].
\end{align*}
We also need the first basic commutator of weight 6: $[d,c,c,c,c,c]$. All the
other basic commutators of weight 6 in $c,d$ lie in $N$. We call these basic
commutators (including $[d,c,c,c,c,c]$) $d_{1},d_{2},\ldots,d_{15}$. Then%
\[
\lbrack d^{r},c^{s}]=d_{3}^{n_{3}}d_{4}^{n_{4}}\ldots d_{15}^{n_{15}}\text{
modulo }N
\]
where $n_{3},n_{4},\ldots,n_{15}$ equal%
\begin{align*}
& rs,r\binom{s}{2},\binom{r}{2}s,r\binom{s}{3},\binom{r}{2}\binom{s}{2}%
,\binom{r}{3}s,r\binom{s}{4},\binom{r}{2}\binom{s}{3},\binom{r}{3}\binom{s}%
{2},\binom{r}{4}s,\\
& 3\binom{r}{2}\binom{s}{3}+2\binom{r}{2}\binom{s}{2}+r\binom{s}{3},4\binom
{r}{3}\binom{s}{2}+2\binom{r}{3}s+3\binom{r}{2}\binom{s}{2}+\binom{r}%
{2}s,r\binom{s}{5}.
\end{align*}
The derivation of these exponents is straightforward, but tedious, so I will
omit it. It is easy to use a computer to check that they are correct for any
number of $r,s$. Using this expression for $[d^{r},c^{s}]$ with $c^{s}%
=[b,a]^{q}$ it is straightforward to check that $[b,a]^{q}$ commutes with all
the elements in $C$ modulo $N$.

This completes our proof that $K/N$ is a direct product of 1375 copies of
the cyclic group of order 8.

We have shown that every element $k\in K$ can be expressed uniquely modulo
$\gamma_{14}(F)$ in the form (6) above. But there is a problem in that the
coefficients $m_3,m_4,\ldots ,m_{1377}$ which appear in (6) 
can depend on $q$. To illustrate this, consider the
following example. Let $c_i,c_j,c_k$ be the basic commutators 
\[
[b,a,a,a,a],\,[b,a,a,a,b],\,[[b,a,a,a,b],[b,a,a,a,a]].
\]
Then working modulo $\gamma_{14}(F)$ we have
\[
(c_ic_j)^{\frac{q}{4}}=c_i^{\frac{q}{4}}c_j^{\frac{q}{4}}c_k^{\frac{q}{8}(\frac{q}{4}-1)}.%
\]
However $8|\frac{q}{4}$ provided $q\geq 32$, and so the representative vector of
$(c_ic_j)^{\frac{q}{4}}N$, thought of as an element in $C_8^{1375}$, is
\[
[0,\ldots,0,1,0,\ldots,0,1,0,\ldots,0,-1,0,\ldots,0]
\]
which does not depend on $q$.

To solve this problem in generality we need to investigate the binomial
coefficients $\binom{q}{d}$ for $1\leq d\leq 13$. These are all divisible 
by $\frac{q}{8}$, so we can write $\binom{q}{d}=\frac{q}{8}n$ for some integer $n$. 
We show that $n\operatorname{mod}8$ only depends on $d$, and not on $q$ (provided
$q\geq32$). 

Consider $\binom{q}{8}$ for example. We want to show that $\binom{q}{8}%
=\frac{q}{8}n$ where $n\operatorname{mod}8$ does not depend on $q$.%
\begin{align*}
\frac{q}{8}n  & =\frac{\allowbreak q^{8}-28q^{7}+322q^{6}-1960q^{5}%
+6769q^{4}-13\,132\allowbreak q^{3}+13\,068q^{2}-5040q}{8!}\\
& =-\frac{q}{8}+\frac{2^{2}\times3267\times q^{2}+rq^{3}}{2^{7}\times315}%
\end{align*}
for some integer $r$. The fraction $\frac{2^{2}\times3267\times q^{2}+sq^{3}%
}{2^{7}\times315}$ is actually an integer, and since $q=2^{k}$ with $k\geq 5 $
we can write this integer as $qm$ for some integer $m$. Dividing through by
$\frac{q}{8}$ we have $n=-1+8m$, and so $n=-1\operatorname{mod}8$.

For another example consider $\binom{q}{12}$, and again write this binomial
coefficient as $\frac{q}{8}n$.%
\[
\frac{q}{8}n=-\frac{q}{3\times4}+\frac{2^{7}\times1024785\times q^{2}+rq^{3}%
}{12!}
\]
for some integer $r$. Multiplying both sides of this equation by 3 we have%
\[
\frac{q}{8}3n=-\frac{q}{4}+sq
\]
for some integer $s$. So $3n\operatorname{mod}8=-2$ and $n\operatorname{mod}8=2$.%

The binomial coefficients $\binom{q}{d}$ where $d$ is odd are all divisible by
$q$, so can all be written in the form $\frac{q}{8}n$ where
$n\operatorname{mod}8=0$. The binomial coefficients $\binom{q}{d}$ for
$d=2,4,6,10$ can be handled in the same way as we dealt with the cases
$d=8,12$.

Similarly the binomial coefficients $\binom{q}{d}$ for $1\leq d\leq7$ are all
divisible by $\frac{q}{4}$ and can all be written in the form $\frac{q}{4}n$
where $n\operatorname{mod}8$ does not depend on $q$.

And finally the binomial coefficients $\binom{q}{d}$ for $1\leq d\leq 3$ are all
divisible by $\frac{q}{2}$ and can all be written in the form $\frac{q}{2}n$
where $n\operatorname{mod}8$ does not depend on $q$.

Now we return to the issue of representing an element $[y^q,x]$ in the form (6)
above. We need to show that provided $q\ge 32$ then the representative vector
of $[y^q,x]N$ depends only on $x,y$, and not on $q$.
We introduce the notation rv$\,(kN)$ for the representative vector of $kN$ when $k\in K$.

Equation (3) from Section 3 gives
\[
\lbrack y^q,x]=[y,x]^{q}[y,x,y]^{\binom{q}{2}}[y,x,y,y]^{\binom{q}{3}}(c_{5}%
\alpha)^{f_{5}(q)}\ldots(c_{1377}\alpha)^{f_{1377}(q)}
\]
modulo $N$, where $\alpha$ is the endomorphism of $F$ mapping $a,b$ to $y,[y,x]$.
And so
\[
\text{rv}\,([y^q,x]N)=\text{rv}\,([y,x]^qN)+\text{rv}\,([y,x,y]^{\binom{q}{2}}N)+
\ldots +\text{rv}\,((c_{1377}\alpha)^{f_{1377}(q)}N).
\]
We show that all the summands on the right hand side of this equation depend only
on $x,y$, and not on $q$.

First consider a summand rv$\,((c_{i}\alpha)^{f_{i}(q)}N)$ where wt$\,c_i\ge 8$.
From our analysis of binomial coefficients we know that $f_i(q)=n\frac{q}{8}$ for
some integer $n$ where $n \bmod 8$ is independent of $q$ (provided $q\ge 32$).
Furthermore $(c_{i}\alpha)\in \gamma_8(F)$ so that $(c_{i}\alpha)^q\in N$. So
\[
(c_{i}\alpha)^{f_{i}(q)}N=(c_{i}\alpha)^{(n\bmod 8)\frac{q}{8}}N.
\]
We show that rv$\,((c_{i}\alpha)^{f_{i}(q)}N)$ depends only on $x,y$ by showing that
if $g\in \gamma_8(F)$ then rv$\,(g^{\frac{q}{8}}N)$ depends only on $g$ and not on $q$. Let
\[
g=c_{42}^{\beta _{42}}\ldots c_{1377}^{\beta _{1377}}\text{ modulo }N,
\]
for some integers $\beta_i$. Then since $\gamma_{8}(F)$ is abelian modulo $N$,
\[
g^{\frac{q}{8}}=c_{42}^{\frac{q}{8}\beta _{42}}\ldots c_{1377}^{\frac{q}{8}%
\beta _{1377}}\text{ modulo }N,
\]
and rv$\,(g^{\frac{q}{8}}N)$ equals
\[
[0,\ldots ,0,\beta _{42},\ldots ,\beta _{1377}]
\]
which depends on $g$, and not on $q$.

Next consider a summand rv$\,((c_{i}\alpha)^{f_{i}(q)}N)$ where $4\le \text{wt}\,c_i<8$.
We can write $f_i(q)=n\frac{q}{4}$ for some integer $n$ where $n \bmod 8$ is independent
of $q$. The element $c_{i}\alpha \in \gamma_5(F)$ and so $(c_{i}\alpha)^{2q}\in N$.
So
\[
(c_{i}\alpha)^{f_{i}(q)}N=(c_{i}\alpha)^{(n\bmod 8)\frac{q}{4}}N.
\]
We show that rv$\,((c_{i}\alpha)^{f_{i}(q)}N)$ depends only on $x,y$ by showing that
if $g\in \gamma_5(F)$ then rv$\,(g^{\frac{q}{4}}N)$ depends only on $g$ and not on $q$. Let
\[
g=c_{9}^{\beta _{9}}\ldots c_{1377}^{\beta _{1377}}\text{ modulo }N,
\]
for some integers $\beta_i$.
Then
\[
g^{\frac{q}{4}}=c_{9}^{\frac{q}{4}\beta _{9}}\ldots c_{1377}^{\frac{q}{4}%
\beta _{1377}}h^{\binom{\frac{q}{4}}{2}}\text{ modulo }N,
\]
where
\[
h=\prod _{9\leq i<j\leq 1377}[c_{j}^{\beta _{j}},c_{i}^{\beta _{i}}].
\]
So
\[
\text{rv}\,(g^{\frac{q}{4}}N)=[0,\ldots ,0,\beta _{9},\ldots ,\beta _{41},%
2\beta _{42},\ldots ,2\beta _{1377}]+u
\]
where $u=\,$rv$\,(h^{\binom{\frac{q}{4}}{2}}N)$. As we have seen
\[
h^{\binom{\frac{q}{4}}{2}}N=h^{-\frac{q}{8}}N
\]
so $u$ depends only on $h$, which in turn depends only on $g$. So
rv$\,(g^{\frac{q}{4}}N)$ depends only on $g$ and not on $q$.

Now consider a summand rv$\,((c_{i}\alpha)^{f_{i}(q)}N)$ where wt$\,c_i=2$ or $3$.
We can write $f_i(q)=n\frac{q}{2}$ for some integer $n$ where $n \bmod 8$ is independent
of $q$. The element $c_{i}\alpha \in \gamma_3(F)$ and so $(c_{i}\alpha)^{4q}\in N$.
So
\[
(c_{i}\alpha)^{f_{i}(q)}N=(c_{i}\alpha)^{(n\bmod 8)\frac{q}{2}}N.
\]
We show that rv$\,((c_{i}\alpha)^{f_{i}(q)}N)$ depends only on $x,y$ by showing that
if $g\in \gamma_3(F)$ then rv$\,(g^{\frac{q}{2}}N)$ depends only on $g$ and not on $q$. Let
$g=c_{4}^{\beta _{4}}h$ modulo $N$ where
\[
h=c_{5}^{\beta _{5}}\ldots c_{1377}^{\beta _{1377}}.
\]
Then from equation (2) in Section 3 we see that
\[
g^{\frac{q}{2}}=c_4^{\frac{q}{2}\beta _{4}}h^{\frac{q}{2}}[h,c_{4}^{\beta _{4}}]^{\binom{\frac{q}{2}}{2}}(c_{4}\gamma)^{f_{4}(\frac{q}{2})}\ldots (c_{1377}\gamma)^{f_{1377}(\frac{q}{2})}\text{ modulo }N
\]
where $\gamma$ is the endomorphism of $F$ mapping $a,b$ to $c_{4}^{\beta _{4}},h$.
If wt$\,c_i>4$ then $c_i\gamma \in \gamma_{14}(F)$ and so
\[
g^{\frac{q}{2}}=c_4^{\frac{q}{2}\beta _{4}}h^{\frac{q}{2}}[h,c_{4}^{\beta _{4}}]^{\binom{\frac{q}{2}}{2}}(c_{4}\gamma)^{f_{4}%
(\frac{q}{2})}\ldots (c_{8}\gamma)^{f_{8}(\frac{q}{2})}\text{ modulo }N
\]
and
\[
 \text{rv}\,(g^{\frac{q}{2}}N)=\text{rv}\,(c_4^{\frac{q}{2}\alpha _{4}}N)+%
 \text{rv}\,(h^{\frac{q}{2}}N)+\ldots +\text{rv}\,((c_{8}\gamma)^{f_{8}(\frac{q}{2})}N).
\]
Clearly rv$\,(c_4^{\frac{q}{2}\alpha _{4}}N)=[0,\alpha _4,0,\ldots,0]$ depends only on $g$.
And we can assume by induction on the length of the product 
$h=c_{5}^{\beta _{5}}\ldots c_{1377}^{\beta _{1377}}$ that rv$\,(h^{\frac{q}{2}}N)$ depends 
only on $h$, and hence only on $g$.

So consider the element $[h,c_{4}^{\beta _{4}}]$. This element lies in $\gamma_6(F)$ so that 
$[h,c_{4}^{\beta _{4}}]^{2q}\in N$. Furthermore
\[
\binom{\frac{q}{2}}{2}=\frac{q}{4}(\frac{q}{2}-1)=-\frac{q}{4}\text{ modulo }2q
\]
provided $q\ge 32$. So
\[
[h,c_{4}^{\beta _{4}}]^{\binom{\frac{q}{2}}{2}}N=[h,c_{4}^{\beta _{4}}]^{-\frac{q}{4}}N
\]
and rv$\,([h,c_{4}^{\beta _{4}}]^{\binom{\frac{q}{2}}{2}}N)$ (as we have seen above) 
depends only on $[h,c_{4}^{\beta _{4}}]$, which in turn depends only on $g$.

For $i=4,5,6,7,8$ $c_i\gamma \in \gamma_9(F)$, so that 
$(c_i \gamma )^q \in N$. It is straightforward to see that $\frac{q}{8}|f_i(\frac{q}{2})$
for $i=4,5,6,7,8$, and it is also straightforward to see that $f_i(\frac{q}{2})$ modulo
$q$ equals $4\frac{q}{8}$, $6\frac{q}{8}$, $7\frac{q}{8}$, $5\frac{q}{8}$, $5\frac{q}{8}$
for $i=4,5,6,7,8$ provided $q\ge 32$. It follows that 
$(c_{4}\gamma)^{f_{4}(\frac{q}{2})}N=(c_{4}\gamma)^{4\frac{q}{8}}N$ and
rv$\,((c_{4}\gamma)^{f_{4}(\frac{q}{2})}N)$ depends only on $c_{4}\gamma$ and hence only on $g$.
Similarly rv$\,((c_{i}\gamma)^{f_{i}(\frac{q}{2})}N)$ depends only on $g$ for $i=5,6,7,8$.

So rv$\,(g^{\frac{q}{2}}N)$ is a sum of vectors each of which depends only on $g$.

Finally consider the summand rv$\,([y,x]^qN)$. Let $[y,x]=[b,a]^{\beta}h$ where $h\in \gamma_3(F)$. 
Then from equation (2) in Section 3 we see that
\[
[y,x]^q=[b,a]^{q\beta}h^q[h,[b,a]^{\beta}]^{\binom{q}{2}}(c_{4}\gamma)^{f_{4}(q)}\ldots %
(c_{23}\gamma)^{f_{23}(q)}\text{ modulo }N.
\]
where $\gamma$ is the endomorphism of $F$ sending $a,b$ to $[b,a]^{\beta},h$. (For 
$i>23$ $c_i\gamma \in \gamma_{14}(F)$). So rv$\,([y,x]^qN)$ equals
\[
[\beta,0,\ldots,0]+\text{rv}\,(h^qN)+\text{rv}\,([h,[b,a]^{\beta}]^{\binom{q}{2}}N)+%
\ldots +\text{rv}\,((c_{23}\gamma)^{f_{23}(q)}N).
\]
We have already shown that rv$\,(h^qN)$ depends only on $h$, and hence only on $x,y$, so consider
rv$\,([h,[b,a]^{\beta}]^{\binom{q}{2}}N)$. Since $[h,[b,a]^{\beta}]\in \gamma_5(F)$
it follows that $[h,[b,a]^{\beta}]^{2q}\in N$. So, as we have seen above,
$[h,[b,a]^{\beta}]^{\binom{q}{2}}N=[h,[b,a]^{\beta}]^{-\frac{q}{2}}N$, so that
rv$\,([h,[b,a]^{\beta}]^{\binom{q}{2}}N)$ depends only on $x,y$. Both $c_{4}\gamma$
and $c_{5}\gamma$ lie in $\gamma _5(F)$, and so $c_{4}\gamma ^{2q}\in N$ and
$c_{5}\gamma ^{2q}\in N$. Now $f_4(q)=\binom{q}{3}$ and $f_5(q)=\binom{q}{2}+2\binom{q}{3}$,
and so $(c_{4}\gamma)^{f_{4}(q)}N=(c_{4}\gamma)^qN$ and 
$(c_{5}\gamma)^{f_{5}(q)}N=(c_{5}\gamma)^{-\frac{q}{2}}N$, and
rv$\,((c_{i}\gamma)^{f_{i}(q)}N)$ depends only on $x,y$ for $i=4,5$.

The elements $c_{i}\gamma$ for $i=6,7,\ldots, 23$ lie in $\gamma_8(F)$, and so all have order
dividing $q$ modulo $N$. And the exponents $f_i(q)$ for $i=6,7,\ldots,23$ can all be
expressed in the form $\frac{q}{8}n_i$ modulo $q$, where $n_i\bmod 8$ does not depend
on $q$ (provided $q\geq 32$). So rv$\,((c_{i}\gamma)^{f_{i}(q)}N)$ depends only
on $x,y$ for $i=6,7,\ldots,23$.

Putting all this together we see that rv$\,([y^q,x]N)$ depends only on $x,y$, and not on $q$.

As we stated earlier in this section, $H=F/M$ where
\[
M=\langle\lbrack y^{q},x]\,:\,x,y\in F\rangle N. 
\]
As $x,y$ range over $F$ the elements rv$\,([y^q,x]N)$ generate an (additive) subgroup $S\le C_8^{1375}$.
An element $k\in K$ lies in $M$ if and only if rv$\,(kN)\in S$. The key point is that provided
$q\geq 32$ the subgroup $S$ does not depend on $q$.

\bigskip
Now consider the claim that $[b,a]$ has order $8q$ in $H$. As stated above, I
have used the $p$-quotient algorithm to confirm this for $q=8,16,32$. For
$q\geq32$ this is equivalent to showing that in $F$, $[b,a]^{8q}\in M$ and
$[b,a]^{4q}\notin M$. We have shown above that $[b,a]^{8q}\in N$. On the other
hand, if $q\geq32$ then $[b,a]^{4q}$ has representative vector $[4,0,0,\ldots
,0]$, and since my computations show that $[b,a]^{4q}\notin M$ when $q=32$
this implies that $[4,0,0,\ldots,0]\notin S$, and hence that $[b,a]^{4q}\notin
M$ for any $q\geq32$.

Next consider the claim that $b^{-q}a^{-q}(ab)^{q}$ has order $8$ in $H$. This
is equivalent to showing that $b^{-8q}a^{-8q}(ab)^{8q}\in M$ and that
$b^{-4q}a^{-4q}(ab)^{4q}\notin M$. From equation (2) in Section 3 we see that%
\[
b^{-8q}a^{-8q}(ab)^{8q}=[b,a]^{\binom{8q}{2}}[b,a,a]^{\binom{8q}{3}}%
c_{5}^{f_{5}(8q)}\ldots c_{1377}^{f_{1377}(8q)}\text{ modulo }N,
\]
and the properties of the integer-valued polynomials $f_{i}(t)$ imply that
\[
b^{-8q}a^{-8q}(ab)^{8q}\in K.
\]
My computer calculations show that
$b^{-8q}a^{-8q}(ab)^{8q}\in M$ for $q=32$. So the representative vector of
$b^{-8q}a^{-8q}(ab)^{8q}$ lies in $S$ when $q=32$. It is not really relevant,
but the representative vector of $b^{-8q}a^{-8q}(ab)^{8q}$ is
\[
[4,0,0,4,4,4,0,0,\ldots,0].
\]
So $b^{-8q}a^{-8q}(ab)^{8q}\in M$ for all
$q\geq32$. The element $b^{-4q}a^{-4q}(ab)^{4q}$ also lies in $K$, and my
computer calculations show that if $q=32$ then $b^{-4q}a^{-4q}(ab)^{4q}\notin
M$. So the representative vector of $b^{-4q}a^{-4q}(ab)^{4q}$ does not lie in
$S$ when $q=32$, and this implies that it does not lie in $S$ for any
$q\geq32$. So $b^{-4q}a^{-4q}(ab)^{4q}\notin M$ for $q\geq32$.

The claims in Theorem 5 and Theorem 7 for the orders of $b^{-q}a^{-q}(ab)^{q}$
and $[b,a]$ in Schur covers of $R(2,q;c)$ for $c<12$ follow similarly. We
replace $N$ by $N\gamma_{c+2}(F)$. If $c_{r}$ is the last basic commutator of
weight $c+1$ then any element in $K/N$ has a unique representative vector
$[m_{3},m_{4},\ldots,m_{r}]$ for all $q\geq32$, and the proof goes through in
the same way as above.

There is a slight problem in showing that $b^{-q}a^{-q}(ab)^{q}$ is
non-trivial in the Schur cover of $R(2,q;c)$ since $b^{-q}a^{-q}(ab)^{q}\notin
K$. But we can directly calculate a Schur cover of $R(2,q;1)$ by hand, and
show that $b^{-q}a^{-q}(ab)^{q}\neq1$ in this cover. Similarly we can show
that $[b,a]^{\frac{q}{2}}$ is non-trivial in a Schur cover of $R(2,q;1)$, so
that the order of $[b,a]$ in a Schur cover of $R(2,q;c)$ is at least $q$ for
any $c$.

Finally, let $G$ be a finite 2-group with exponent $q>2$ and nilpotency
class at most 3. We show that if $H$ is the covering group of $G$ then
$H'$ has exponent dividing $q$. Our calculations with the covering group
of $R(2,q;3)$ show that $[y,x]^q=[y,x,x]^{\frac{q}{2}}=1$ for all $x,y\in H$.
Groups satisfying the 2-Engel identity $[y,x,x]=1$ are nilpotent of class
at most 3. So if $h\in \gamma _4(H)$, $h$ can be expressed as a product
of terms $[y,x,x]$ and their inverses, with $x,y\in H$. So $h^{\frac{q}{2}}$
is a product of terms $[y,x,x]^{\frac{q}{2}}$ and their inverses, and
hence $h^{\frac{q}{2}} = 1$. This implies that $\gamma_4(H)$ has exponent 
dividing $\frac{q}{2}$. Since $H'$
is generated by commutators which have order dividing $q$,  and
$\gamma_4(H)$ has exponent dividing $\frac{q}{2}$, we see that
that $H'$ has exponent dividing $q$.

\section{Proofs of Theorem 6 and Theorem 8}

The proofs of Theorem 6 and Theorem 8 are essentially the same as the proofs
of Theorem 5 and Theorem 7. Let $q=3^{k}$ where $k\geq2$, let $R(2,q;12)=F/R $
where $F$ is the free group of rank two generated by $a,b$, and let
$H=F/[F,R]$. We can use the $p$-quotient algorithm to show that Theorem 6 and
Theorem 8 hold true when $q=9$ or 27, and we show that the fact that they hold
true for $q=27$ shows that they also hold true for higher powers of 3.

We let $q=3^{k}$ where $k\geq3$, and let $c_{1},c_{2},\ldots,c_{1377}$ be our
list of basic commutators of weight at most $13$. As in the proof of Theorem 5
and Theorem 7 we let $x,y\in H$ and obtain a relation%
\begin{equation}
\lbrack y,x]^{q}[y,x,y]^{\binom{q}{2}}[y,x,y,y]^{\binom{q}{3}}(c_{5}%
\alpha)^{f_{5}(q)}\ldots(c_{1377}\alpha)^{f_{1377}(q)}=1
\end{equation}
where $\alpha$ is the homomorphism from $F$ to $H$ mapping $a,b$ to $y,[y,x]$. 
If wt$\,c_{i}=w$ then $c_{i}\alpha$ is a commutator in $x$ and $y $, with
$w$ entries $y$, and $\deg f_{i}\leq w$. The binomial coefficients $q$
and $\binom{q}{2}$ are both divisible by $q$, the binomial coefficients
$\binom{q}{d}$ for $d<9$ are divisible by $\frac{q}{3}$, and the
binomial coefficients $\binom{q}{d}$ for $d<27$ are divisible by
$\frac{q}{9}$.

If we let $y\in\gamma_{5}(H)$ then all commutators in $H$ with 3 or more entries $y$ 
are trivial and we obtain the relation $[y,x]^{q}[y,x,y]^{\binom{q}{2}%
}=1$. This implies that $[y,x]^{q}=1$, and so (since $\gamma_{6}(H)$ is nilpotent
of class 2) $\gamma_{6}(H)$ has exponent $q$.

Now let $y\in\gamma_{2}(H)$, and replace $q$ by $3q$ in equation (7). Using
the fact that $\gamma_{6}(H)$ has exponent $q$ we see that $[y,x]^{3q}%
[y,x,y]^{\binom{3q}{2}}=1$, which implies that $[y,x]^{3q}=1$, and hence that
$\gamma_{3}(H)$ has exponent $3q$.

Finally replace $q$ by $9q$ in equation (7) and using the facts that
$\gamma_{3}(H)$ has exponent $3q$ and $\gamma_{6}(H)$ has exponent $q$ we see 
that $[y,x]^{9q}=1$ for all $x,y\in H$, and that $H^{\prime}$ has exponent $9q$.

Let $N$ be the normal subgroup 
\[
\gamma_{2}(F)^{9q}\gamma_{3}(F)^{3q}\gamma_{7}(F)^{q}\gamma_{14}(F)<F.
\]
Then $H=F/M$ where $M=\langle\lbrack y^{q},x]\,:\,x,y\in F\rangle N$.

Now let
\[
K=\gamma_{2}(F)^{q}\gamma_{3}(F)^{\frac{q}{3}}\gamma_{7}(F)^{\frac{q}{9}}\gamma_{14}(F).
\]
Just as in Section 5 we can show that every element $k\in K$ can be uniquely expressed
modulo $\gamma_{14}(F)$ in the form
\begin{equation}
k=[b,a]^{qm_{3}}[b,a,a]^{\frac{q}{3}m_{4}}\ldots c_{23}^{\frac{q}{3}m_{23}%
}c_{24}^{\frac{q}{9}m_{24}}\ldots c_{1377}^{\frac{q}{9}m_{1377}}%
\end{equation}
(There are 23 basic commutators $c_{i}$ of weight at most 6.) Similarly every element
of $N$ can be uniquely expressed modulo $\gamma_{14}(F)$ in the form
\[
[b,a]^{9qn_{3}}[b,a,a]^{3qn_{4}}\ldots c_{23}^{3qn_{23}%
}c_{24}^{qn_{24}}\ldots c_{1377}^{qn_{1377}}.%
\]

If $k\in K$ is given by equation (8) then $k\in N$ if and only if $9|m_i$ for 
$i=3,4,\ldots ,1377$. And just as in Section 5 we can show that $K/N$ is abelian
and is a direct product of 1375 copies of the cyclic group of order 9.
We let $[m_{3},m_{4},\ldots,m_{1377}]$ be the representative
vector for $kN$, and think of this vector as an element in $C_{9}^{1375}$. 
Multiplying elements of $K/N$ corresponds to adding their representative
vectors, and the element $k$ lies in $N$ if and only if the representative 
vector of $kN$ is $0$.

A similar argument to that given in Section 5 for binomial coefficients
$\binom{q}{d}$ $(d\leq13)$ for $q$ a power of 2 at least as big as 32, shows
that if $q\geq27$ is a power of 3 then all the binomial coefficients
$\binom{q}{d}$ $d\leq13$ can be expressed in the form $\frac{q}{9}n$ for some
integer $n$ where $n\operatorname{mod}9$ depends only on $d$, and not on $q$.
Similarly, if $q\geq27$ is a power of 3 then all the binomial coefficients
$\binom{q}{d}$ $d<9$ can be expressed in the form $\frac{q}{3}n$ for some
integer $n$ where $n\operatorname{mod}9$ depends only on $d$, and not on $q$.
And finally, if $q\geq27$ is a power of 3 then 
$\binom{q}{1}=q$ and $\binom{q}{2}=qn$ where $n \bmod 9=4$.
Using the same argument as we used in Section 5, we see that if $q\geq27$ and
$x,y\in F$, then $[y^{q},x]\in K$, and the representative vector of
$[y^{q},x]N$ depends only on $x,y$, and not on $q$. As we stated above,
$H=F/M$ where $M=\langle\lbrack y^{q},x]\,:\,x,y\in F\rangle N$. The quotient
$M/N$ is a subgroup of $K/N$ and
the set of representative vectors of elements in this subgroup is a subgroup
$S\leq C_{9}^{1375}$. If $k\in K$, then $k\in M$ if and only if the
representative vector of $kN$ lies in $S$.

The remainder of the proof of Theorem 6 goes through in the same
way as in Section 5, as does the proof of the claims in Theorem 8
for the order of $[b,a]$ in Schur covers of $R(2,q;c)$ for various $c$.

Now let $G$ be a finite 3-group with exponent $q>3$ and class at most 8,
and let $H$ be the covering group for $G$. We show that $H'$ has exponent
dividing $q$. We have shown that commutators in $H$ have order dividing
$q$, but we need to show that products of commutators have order
dividing $q$. Our calculations with the covering group of $R(2,q;8)$
show that $[y,x,x,x]^{\frac{q}{3}}=1$ for all $x,y\in H$. It is known
that 3-Engel groups are locally nilpotent, and Werner Nickel's nilpotent
quotient algorithm in \textsc{Magma} \cite{boscan95} has a facility
for computing Engel groups. The free 3-Engel group of rank 5 has class 9,
and it is an easy calculation with the nilpotent quotient algorithm
to show that 3-Engel groups satisfy the identity $[x_1,x_2,x_3,x_4,x_5]^{20}=1$.
This implies that in a free group $[x_1,x_2,x_3,x_4,x_5]^{20}$ can be expressed
as a product of terms $[y,x,x,x]$ and their inverses. Now $\gamma_5(H)$
is abelian and is generated by elements with order dividing $q$. So
$\gamma_5(H)$ has exponent dividing $q$, which is coprime to 20.
So if $h\in \gamma_{5}(H)$ then $h$ can be expressed as a product
of terms $[y,x,x,x]$ and their inverses (with $x,y\in H$).
This implies that $h^{\frac{q}{3}}$ is a product of terms $[y,x,x,x]^{\pm \frac{q}{3}}$
(which are all trivial) and terms
\[
[[y,x,x,x]^{\pm 1},[z,t,t,t]^{\pm 1}]^{\binom{\frac{q}{3}}{2}}
\]
which are also trivial. So $\gamma_5(H)$ has exponent dividing $\frac{q}{3}$.
This, combined with the fact that $H'$ is generated by elements with
order dividing $q$, implies that $H'$ has exponent dividing $q$.

\section{Proof of Theorem 9 and Theorem 10}

To prove Theorem 9 we need to show that if $q$ is a power of $5$ then the order of 
$[b,a]$ in a Schur cover of $R(2,q;c)$ is $q$ for $c<9$, and $5q$ for $c=9$. It is an 
easy calculation with the $p$-quotient algorithm to show that this is the case for
$q=5,25$. We use the same argument as in Section 5 and Section 6 to show that
this implies that the theorem holds true for all powers of 5.

Let $q=5^{k}$ where $k\geq 2$, let $R(2,q;9)=F/R$ where $F$ is the free group
of rank two generated by $a,b$, and let $H=F/[F,R]$. Let $c_{1},c_{2}%
,\ldots,c_{226}$ be our list of basic commutators of weight at most $10$. As
in the last two sections we let $x,y\in H$ and obtain a relation%
\begin{equation}
\lbrack y,x]^{q}[y,x,y]^{\binom{q}{2}}[y,x,y,y]^{\binom{q}{3}}(c_{5}%
\alpha)^{f_{5}(q)}\ldots(c_{226}\alpha)^{f_{226}(q)}=1
\end{equation}
where $\alpha$ is the homomorphism from $F$ to $H$ mapping $a,b$ to $y,[y,x]
$. If wt$\,c_{i}=w$ then $c_{i}\alpha$ is a commutator in $x$ and $y$, with
$w$ entries $y$, and $\deg f_{i}\leq w$. The binomial coefficients
$\binom{q}{d}$ for $d<5$ are divisible by $q$, and the binomial
coefficients $\binom{q}{d}$ for $d<25$ are divisible by $\frac{q}{5}$.

If we let $y\in\gamma_{2}(H)$ then all commutators in $H$ with 5 or more entries $y$
are trivial and we see that $[y,x]^{q}$ is a product of $q^{th}$ powers
of commutators of higher weight in $x,y$. So $[y,x]^{q}=1$, and
$\gamma_{3}(H)$ has exponent $q$.

Now replace $q$ by $5q$ in (9) and we obtain $[y,x]^{5q}=1$, which implies
that $\gamma_{2}(H)$ has exponent $5q$.

So we let
\[
N=\gamma_{2}(F)^{5q}\gamma_{3}(F)^{q}\gamma_{11}(F)
\]
and we let
\[
K=\gamma_{2}(F)^{q}\gamma_{3}(F)^{\frac{q}{5}}\gamma_{11}(F).
\]
Just as in Section 5 and Section 6 we can show that $K/N$ is a direct
product of 224 copies of the cyclic group of order 5. Every element $k\in K$
can be uniquely expressed modulo $N$ in the form
\[
k=\lbrack b,a]^{qm_{3}}c_{4}^{\frac{q}{5}m_{4}}\ldots c_{226}^{\frac{q}{5}m_{226}}%
\]
where $0\leq m_{i}<5$ for $i=3,4,\ldots,226$. We let $[m_3,m_4,\ldots ,m_{226}]$
be the representative vector for $kN$, and we think of it as an element
in $C_5^{224}$.

Just as in Section 5 we can show that if $x,y\in F$ then $[x^q,y]\in K$, and the
representative vector of $[x^q,y]N$ depends only on $x,y$, and not on $q$.
The rest of the proof that the order of $[b,a]$ in a Schur
cover of $R(2,q;c)$ is $q$ for $c<9$, and $5q$ for $c=9$ goes through
just as in Section 5.

If $G$ is any group of exponent $q=5^k$ ($k\ge 1$) with class less than 9, and if $H$ is the
cover of $G$ then commutators in $H$ have order dividing $q$, and so (since
$H'$ has class at most 4) $H'$ has exponent dividing $q$, which implies
that $M(G)$ has exponent dividing $q$.

The proof of Theorem 10 is almost identical to the proof of Theorem 9.

\end{document}